\documentclass[10pt]{article}
\usepackage{latexsym,amssymb,verbatim}
\usepackage{enumerate}
\usepackage{graphics}

\usepackage{latexsym,amssymb}
\usepackage{enumerate, verbatim, mathrsfs}
\usepackage{graphics}

\usepackage{amsthm}
\usepackage{amssymb}
\usepackage{latexsym}
\usepackage{longtable}
\usepackage{epsfig}
\usepackage{amsmath}
\usepackage{hhline}

\newtheorem{bthm}{\sc Theorem}
\newtheorem{thm}{\sc Theorem}[section]
\newtheorem{cor}[thm]{\sc Corollary}
\newtheorem{lem}[thm]{\sc Lemma}
\newtheorem{prop}[thm]{\sc Proposition}
\newcommand{\pf}{{\it Proof:\quad}}
\newcommand{\dne}{\hfill $\Box$\smallskip}
\newcommand{\FF}{\mathbb{F}}

\newcommand{\CP}{{\cal P}}

\newcommand{\II}{{\rm I}}

\newcommand{\p}{\partial}

\newcommand{\set}[1]{\left\{#1\right\}}
\newcommand{\gd}{{\partial}}
\newcommand{\pa}{{\partial}}
\newcommand{\ga}{{\alpha}}

\newcommand{\gb}{{\beta}}
\newcommand{\ic}{{\rm inc}}
\newcommand{\lra}{\,\,\leftrightarrow\,\,}
\newcommand{\ra}{\hookrightarrow}

\begin{document}

\title{Incidence Homology of Finite Projective Spaces}
\author{Johannes Siemons {\it \,and\,\,} Daniel Smith\\
\\ School of Mathematics \\University of East Anglia, Norwich, UK\\
}
\date{\scriptsize Version of February 2012, printed \today}

\maketitle
\begin{abstract}
\noindent
Let $\FF$ be the  finite field of $q$ elements and let $\CP(n,q)$ denote the projective space of dimension $n\!-\!1$ over  $\FF.$  We construct a family $H^{n}_{k,i}$  of combinatorial homology modules associated to  $\CP(n,q)$ for  a coefficient   field  $F$ of positive  characteristic co-prime to $q.$ As $F{\rm GL}(n,q)$-representations  these modules are obtained from the  permutation action of ${\rm GL}(n,q)$ on  the set of  subspaces of $\FF^{n}.$ We prove a branching rule for $H^{n}_{k,i}$  and use this to determine the homology representations completely. Results include a duality theorem, the characterisation of  $H^{n}_{k,i}$ through the standard irreducibles of ${\rm GL}(n,q)$ over $F$ and applications.


\end{abstract}

\section{Incidence Homology}

Let $\CP$ be a finite ranked partially ordered set and let $F$ be a field. In this paper we show that there is  a family of
 homology modules with $F$ as coefficient domain that can be associated to $\CP.$ This homology contains significant combinatorial and algebraic information when $\CP$ is a finite projective space.  We determine these modules  completely for this case when the coefficient field has  positive characteristic  co-prime to the characteristic of the space.   
 
The homology modules appear in the following fashion. We assume that the rank function ${\rm rk}\!:\CP\to\mathbb{N}\cup\{0\}$ is adjusted  so that $\min\{{\rm rk}(x)\,:\,x\in \CP\}=0.$ For the integer $k\geq 0$ let $\CP_{k}$ denote  the elements of rank  $k$ in  $\CP$ and let $M_{k}:=F{\cal P}_{k}$ be  the $F$-vector space with basis ${\cal P}_{k};$ in particular $M_{k}=0$ if $\CP_{k}=\emptyset.$ 

The partial order on $\CP$ now provides a linear {\it incidence map } $\partial\!: M_{k}\to M_{k-1}$ defined  by $\partial(x)=\sum y$ for $x\in \CP_{k}$ where the sum runs over all $y\in \CP_{k-1}$ covered by $x.$ This gives rise to the sequence  
\begin{equation}\label{SEQ001}
{\mathcal M}:\quad 0\stackrel{\,\,\,\partial}{\leftarrow}M_{0}\stackrel{\,\,\,\partial}{\leftarrow}
\ldots\stackrel{\,\,\,\partial}{\leftarrow}
M_{k-2}^{}\stackrel{\,\,\,\partial}{\leftarrow}M_{k-1}\stackrel{\,\,\,\partial}{\leftarrow}M_{k}\stackrel{\,\,\,\partial}{\leftarrow}M_{k+1}
\stackrel{\,\,\,\partial}{\longleftarrow}\ldots \stackrel{\,\,\,\partial}{\longleftarrow} 0\,.
\end{equation}
Since $\CP$ is finite  ${\mathcal M}$ has only finitely many non-zero terms, and in particular, $\pa$ is nilpotent. Therefore we may choose  an integer $m$ such that $\pa^{m}=0$ on $\bigoplus_{k\in\mathbb{Z}}\,M_{k}.$ The options  for $m$ depend on both $\CP$ and $F.$  If we fix  some $k$ and  $0<i<m$  then (\ref{SEQ001}) gives rise to the `subsequence' 
\begin{equation}\label{SEQkii}
{\mathcal M}^{}_{k,i}:\quad 0\stackrel{\,\,\,\partial^*}{\leftarrow}
\ldots\stackrel{\,\,\,\partial^*}{\leftarrow}
M_{k-m}^{}\stackrel{\,\,\,\partial^*}{\leftarrow}
 M_{k-i}^{}\stackrel{\,\,\,\partial^*}{\leftarrow}
M_{k}^{}
\stackrel{\,\,\,\partial^*}{\leftarrow} M^{}_{k+m-i}
\stackrel{\,\,\,\partial^*}{\leftarrow} M^{}_{k+m}
\stackrel{\,\,\,\partial^*}{\leftarrow}\ldots\ \stackrel{\,\,\,\partial^*}{\leftarrow}0
\end{equation}
in which $\partial^{*}$ denotes the appropriate power of $\partial$.  Since  $\partial^{*}\partial^{*}=\pa^{m}=0$ it follows that  ${\mathcal M}_{k,i}$ is homological. For instance, when $m=2$ then ${\mathcal M}={\mathcal M}_{k,i}$ is homological  in the usual way. The
homology at $M_{k-i}^{}\leftarrow M_{k}^{}\leftarrow M^{}_{k+m-i}\,,$ denoted 
by
$$
H_{k,i}^{}:=(\ker\,\partial^{i}\cap M_{k}^{})\,\,\Big/\,\,
                      \partial^{m-i}(M_{k+m-i}^{})\,,
$$
is the {\it  incidence homology } of ${\cal P}$ with coefficient field  $F$ for parameters $k$ and $i.$ This construction is canonical in the sense that if $G$ is a group of automorphisms of $\CP$ then  $H_{k,i}^{}$ is an $FG$-module. This homology has appeared in various guises before, often  under additional restrictions. We mention only some papers~\cite{Mayer, Kapra, fisk, Mie, Mori,  MSP, Fib}. The full details of this construction are explained in Section~3.

\bigskip
In this paper we  are interested in  the homology  of finite projective spaces. Let $q$ be a prime power and  $\FF$ the field of $q$ elements. For  $n\geq 0$ 
let  $\CP=\CP(n,q)$ be the projective space of dimension $n-1$ over $\FF.$  As a partially ordered set  $\CP$ consists of all subspaces of $\FF^{n}$ ordered by containment. 

We assume that the coefficient field $F$ has  positive  characteristic $p_{}$ not dividing  $q.$ Here the natural choice for $m$ is the  quantum characteristic $m=m(p_{},q)$ of $q$ in $F,$ see Section~2. 
We denote $G_{n}:={\rm GL}(n,q)$ and let $S_{n-1}$ be a Singer cycle of $G_{n-1},$ of order $|S_{n-1}|=q^{n-1}-1.$ For $0<i<m=m(p_{},q)$ denote $H_{k,i}^{n}:= H_{k,i}$ and put  $H_{k,0}^{n}=0=H_{k,m}^{0}.$ The following branching rule is familiar from James' book~\cite{GDJ}, it is the key to the incidence  homology.

\medskip
\begin{bthm} [Homology Decomposition]\label{bthm:HMD1}
 \, \, Let $H^{n}_{k,i}$ denote the incidence homology of $\CP(n,q)$ over a field $F$ of characteristic $p_{}>0$ not dividing $q.$  Suppose that $0 \leq  k \leq n$ with $1\leq n$ and $0<i<m=m(p_{},q).$  Then
\[H^n_{k,i}\,\, \cong \,\,H^{n-1}_{k,i+1}\,\, \oplus\,\, H^{n-1}_{k-1,i-1} \,\,\oplus\,\, H^{n-2}_{k-1,i}\,S_{n-1}\]
as $FG_{n-1}$-modules.  \end{bthm} 

We prove a more general version of this theorem in Section~3.  It contains an important rule for the index pairs $(k,i).$ We call $(k,i)$ a {\it middle index for $n$} provided that $n<2k+m-i<n+m.$ This inequality for the parameters is preserved when passing from one side of the isomorphism to the other. By induction we are able to conclude that $H^{n}_{k,i}$ vanishes unless $(k,i)$ is a middle index, see Theorem~\ref{bthm:Comp}\, below.  The main properties of the homology of  projective spaces  in co-prime characteristics are summed up in the following two theorems.  

\medskip
\begin{bthm} [Duality]\label{bthm:Dual} Assume that  $(k,i)$ is a middle index for $n$ and let  $j:=2k-n+m-i$ where $m=m(p,q).$  Then
\vspace{-6mm}
\begin{enumerate}[\!\!\!(i)\,\,]
\item  $H^n_{k,i}\,\, \cong \,\,H^{n}_{n-k,m-i}$ and \vspace{-2mm}
\item $H^n_{k,i}\,\, \cong \,\,H^{n}_{k,j}$
\end{enumerate}

\vspace{-0.5cm}
as $FG_{n}$-modules.  \end{bthm} 

The incidence homology therefore has a $C_{2}\times C_{2}$ symmetry, we give an example to illustrate this in Section~5\, where the theorem is proved. The first part can be interpreted as saying that the duality between $k$-dimensional and $(n-k)$-dimensional subspaces in projective space remains in place for the  homology. At a formal level this can also be understood as a Poincar\'{e} duality. The second duality appears to have no immediate geometric interpretation as far as we are aware.  

In  Section~5\, we provide  the complete decomposition of the $H^{n}_{k,i}$ into standard irreducibles of ${\rm GL}(n,q)$ over $F.$ To state this result let  $(k,i)$ be  a middle index for $n$ and define the following parameter intervals 

\vspace{-2mm}
 {\small $$\begin{array}{rlll} 
T_{k,i}&:=&\{\,t\!:\,k\leq t\leq n-k+i-1\}\,&\text{if $k\geq \frac 12 n$ and  $i\leq \frac 12 (m-n) +k,$}\smallskip\\\smallskip
T_{k,i}&:=&\{\,t\!:\,k\leq t\leq k+m-i-1\}\,&\text{  if $k\geq \frac 12 n$ and  $i>\frac 12 (m-n) +k,$}\\\smallskip
T_{k,i}&:=&\{\,t\!:\,n-k\leq t\leq n-k+i-1\}\,&\text{if $k<\frac 12 n$ \,and\, $i\leq \frac 12 (m-n)+k,$ and }\\\smallskip
T_{k,i}&:=&\{\,t\!:\,n-k\leq t\leq k+m-i-1\}\,&\text{if $k<\frac 12 n$ \,and\, $i>\frac 12 (m-n)+k.$} \end{array}$$}

\vspace{-4mm}
Let $\lambda$ be a composition of $n$ and let  $D^{\lambda}$ denote the head of the Specht module $S^{\lambda}$ over $F.$ In Sections~4.1~and~5.2 we prove 
  
\medskip
\begin{bthm} [Irreducibles]\label{bthm:Comp} Let $0 \leq  k \leq n$ and $0<i<m=m(p_{},q).$  Then
\vspace{-5mm}
\begin{enumerate}[\!\!\!(i)\,\,]
\item $H^n_{k,i}\neq 0 $ if and only if $(k,i)$ is a middle index; \vspace{-2mm}
\item Let $(k,i)$ be a middle index. Then $H^{n}_{k,i}=\bigoplus\, D^{(n-t,t)} $ where the summation runs over all $t$  in $T_{k,i}.$
\end{enumerate}
\end{bthm}

\vspace{-0.2cm}
The first part of this theorem appeared already in~\cite{MSP}, the proof here is more direct. We mention that all three theorems above remain true in the limit case $`q=1\!$'\, when the general linear group ${\rm GL}(n,\FF)$ is replaced by the symmetric group ${\rm Sym}(n)$ and  the 'Singer cycle' of order $q^{n-1}-1=0$ vanishes.  

\bigskip
In the case of finite projective spaces the sequence ${\cal M}_{k,i}$ in \,(\ref{SEQkii})\, has a remarkable property: For every $(k,i)$ it is {\it almost exact,}\, in the sense that it is exact in all but at most one position. A homology arising  from ${\cal M}_{k,i}$  is non-zero precisely when its index parameters satisfy the middle index condition mentioned before; the details are  given in Section 4.

We mention several consequences of this fact.  From a standard application of the trace formula it follows that every non-zero homology appears as an alternating sum (in the Burnside ring) of the permutation modules involved in \,(\ref{SEQkii}), see Corollary~\ref{cor:trace}. This in particular provides an explicit character formula for every irreducible ${\rm GL}(n,q)$-representation appearing in Theorem~\ref{bthm:Comp}\, in terms of permutation characters, and it also gives an explicit formula for Betti numbers. This is shown  in Theorem~\ref{thm:Betti}. 
In the same context we  mention also the rank modulo $p$ of the incidence matrix of $s$- versus $t$-dimensional subspaces of $\FF^{n}.$ This has been determined in Frumkin and Yakir~\cite{Frum}. In a forthcoming paper~\cite{JSV} we show that this question  has a natural interpretation in terms of the incidence homology of  a certain rank-selected poset obtained from the projective space. Other applications~\cite{MS2} concern the multiplicities of irreducible $FG$-characters when $G$ is an arbitrary subgroup of ${\rm GL}(n,\FF)$ acting on subspaces of $\FF^{n}.$ This question was raised by Stanley~\cite{Stan82} for ordinary representations over $\mathbb{C}.$

The methods of the paper are elementary, and we use the standard theory of  ${\rm GL}(n,\FF)$ representations in cross-characteristics from James' book~\cite{GDJ} wherever possible.  We  thank Alex Zalesskii for helpful comments  on an earlier version of this paper. 

\section{Notation and Prerequisites}\label{Sec2}

\vspace{-3mm}
Let $q$ be a prime power and  $\FF$ the field of $q$ elements   and  let  $n\geq 0$ be an integer.  For the integer $i\geq 1$
let $[i]_q:=1+q+\cdots+q^{i-1}.$ Then $(i!)_{q}:=[i]_{q}\cdot [i-1]_{q}\cdot ...\cdot [1]_{q}$ is the $q$-factorial of $i.$ If $n\geq k\geq 0$ are integers then the $q$-binomial coefficient, or Gaussian polynomial, is denoted by 
\[{n \choose k}_{\!q} := \frac{[n]_q \cdot [n-1]_q \cdots 
[n-k+1]_q}{[k]_q \cdot [k-1]_q \cdots [1]_q}=\frac{[(n)!]_{q}}{[(n-k)!]_{q}\cdot[(k)!]_{q}}\,.\] This  is the number of $k$-dimensional subspaces of $\FF^{n}.$

\medskip
Let $p$ be a prime not dividing $q$ and $F:={\rm GF}(p).$ Then the integer 
\begin{equation}\label{quantum}
m=m(p,q):=\min\{\,i>1\,:\, [i]_{q}=0 \,\,{\rm in}\,\,\,F\}
\end{equation}
is the {\it (quantum) characteristic\,} 
 of $\FF$ in $F.$  In other words, $m$ is the  order of $q$ in $F^{\times}$ if $p$ does not divide $q-1$ while   $m=p$ if $p$ divides $q-1.$ In particular, $m$ is the least integer such that $(m!)_{q}=0$ in $F.$

\subsection{Permutation action of ${\rm GL}(n,\FF)$ on the set of subspaces of $\FF^{n}$} \label{Sec2.1}
\vspace{-0.3cm}

Let $n,$ $q$ and $\FF$ be as above and denote $G_{n}:={\rm GL}(n,q).$ Then $G_{n}$ acts on $V:=\FF^{n}$ (written as row vectors) via the matrix product $g\!:v\mapsto vg$ for $v$ in $V$ and  $g$ in  $G_{n}.$ We view $G_{n-1}$ as the subgroup  
\[G_{n-1}=\set{\left({\small\begin{array}{cc}1&0\\0&g' \end{array}}\right)\,\,:\,\, g' \in  {\rm GL}(n-1,q)\,\,}\]
in $G_{n}$ and embed $G_{n-2}\subseteq G_{n-1}$ correspondingly.  The affine linear group $A_{n-1}$ is the subgroup 
$$A_{n-1}:=\set{\left({\small \begin{array}{cc}1&0\\a&g'\end{array}}\right)\,\,:\,\,a \in (\FF^{n-1})^{T},\,\,g' \in G_{n-1}\,\,}\,.$$ For matrices it is our convention that $1$ always stands for  a $1\times 1$ submatrix while a $0$-entry stands for an appropriate column or row of zeros.

\medskip
Let $v_{1},\,v_{2},\,...,\,v_{n}$  be the standard basis of $V.$
If  $x$ is a subspace of $V$ with basis $x_{1},...,x_{k}$ we express the basis vectors  as $x_{i}=\sum_{i=1}^{n}\,x_{ij}\,v_{j}$ in terms of the standard basis. In this way $x$ is represented by the $k\times n$ matrix $(x_{ij})$ of rank $k. $ It is clear now  that $G_{n}$ acts on the set of subspaces of $V$ via the matrix product
$$g\!:\,x\mapsto xg=(x_{ij})g\,\,.$$
Another $k\times n$ matrix $(x'_{ij})$ represents the same $x$ for some other basis  $x'_{1},...,x'_{k}$ if and only if  
$(x'_{ij})$ is of the form  $(x'_{ij})=h\,(x_{ij})$ for some $h$ in ${\rm GL}(k,q).$

For $0\leq k\leq n$ let $L^{n}_{k}$ denote  the set of all subspaces of dimension $k$ in $V=\FF^{n}.$  It is convenient to specify

\vspace{-3mm}
$$\begin{array}{lll}
V^{n-1}&:=\langle v_{2},...,\,v_{n}\rangle\,\,,    \qquad \qquad
&V^{n-2}:=\langle v_{3},...,\,v_{n}\rangle \text{\quad and}\nonumber\bigskip\\
L^{n-1}_{k}&:=\{\,x\in L^{n}_{k}\,:\,x\subseteq V^{n-1}\,\}\,, \qquad
&L^{n-2}_{k}:=\{\,x\in L^{n}_{k}\,:\,x\subseteq V^{n-2}\,\}\,. \nonumber
\end{array} $$

\subsection{Permutation modules of ${\rm GL}(n,q)$ on subspaces in cross characteristic}\label{Sec2.2}

\vspace{-3mm}
As above let $F$ be a field of characteristic ${\rm char}(F)=p_{}>0$ where we assume that  $p_{}$ does not divide $q=|\FF|.$   If $X$ is an arbitrary set we denote by $FX$  the $F$-vector space with basis $X.$ If a group $G$ acts on $X$ then   $FX$  is the $FG$-permutation module afforded by $G$  over $F.$ In our situation $G$ will be  $G_{n}={\rm GL}(n,q),$ or a subgroup of it, and $X$ will be some collection of subspaces of $V=\FF^{n}.$ This section follows James' book~\cite{GDJ}\, closely. We  denote
$$M^{n}_{k}:=FL^{n}_{k}$$  and set $M^{n}_{k}=0$ when $k<0$ or $k>n.$ Similarly let 
\begin{eqnarray}\label{L^{n-1}}
M^{n-1}_{k}&:=FL^{n-1}_{k}\,\,\,\,\,\,\mbox{and} \,\,\,\,\,
M^{n-2}_{k}&:=FL^{n-2}_{k}
\end{eqnarray} with  $M^{n-1}_{k}=0=M^{n-2}_{k}$ when $L^{n-1}_{k}=\emptyset$ or $L^{n-2}_{k}=\emptyset.$ It is clear that $M^{n}_{k}\supseteq M^{n-1}_{k}\supseteq M^{n-2}_{k}$ are modules for $FG_{n},$ $FG_{n-1}$ and $FG_{n-2}$ respectively.

\bigskip
Next we assume that $F^{\times}$ contains an element of order  $q_{0}$ where $q_{0}$ is the prime dividing $q.$  This implies that  there is  a non-trivial homomorphism $$\chi\!:(\FF,+)\rightarrow F^{\times}\,.$$
In  $FG_{n}$ consider the elements 
{$$E_{1}:=\frac{1}{q^{n-1}}\,\,\,\sum_{\ga\in\FF^{n-1}}\,\left({\small \begin{array}{cc}
1 &   0   \\
\ga& \,\II_{n-1}\, \end{array}}\right)$$}
and $$
E_{2}:=\frac{1}{q^{2n-3}}\,\sum_{\ga_{1}\in\FF,\ga,\gb\in\FF^{n-2}}\,\chi(\ga_{1})
\left({\small \begin{array}{ccc}
     1    & 0             &   0     \\
\ga_1 & 1             &   0     \\
\ga &  \gb&   \II_{n-2}       \end{array}}\right)\,\,.
$$
Here $\II_{n-1}$ and $\II_{n-2}$ are identity matrices, and we emphasize  that $E_{1}$ and $E_{2}$ are  elements of $FG_n$ rather than matrix sums.  It is easy to show that $E_{1}E_{2}=0=E_{2}E_{1}$ and $E_iE_i=E_i$ for $i=1,\,2.$

The elements of $A_{n-1}$ commute with $E_1$ and so, if $M$ is an $FG_{n}$-module, then $ME_1:=\langle xE_1 : x \in M \rangle_{F}$ is an $FA_{n-1}$-module. Similarly, 
$ME_2A_{n-1}=\langle xE_2g \,\,: \,\,x \in M,\, g \in A_{n-1} \rangle_{F}$ is an  $FA_{n-1}$-module and this can be written as $ME_2A_{n-1}=ME_{2}S_{n-1}$ when $S_{n-1 }$ denotes a Singer cycle of $G_{n-1}.$  
The following theorem of James is the essential branching rule for modules involved  in $M^{n}_{k}.$

\medskip
\begin{thm} [Theorems 9.11 and 10.5\,~in~\cite{GDJ}] \label{thm:BR}
Assume that $F^{\times}$ contains an element of order  $q_{0}$ where $q_{0}$ is the prime dividing $|\FF|.$
If $M$ is an $FG_n$-module involved in $FL^{n}_{k}$ then\newline
(i) \,\, $M =ME_1 \,\,\oplus \,\, ME_2S_{n-1}$\,\, as $F\!A_{n-1}$-modules, and \newline
(ii) \, $\dim M = \dim(ME_{1})+(q^{n-1}-1) \dim (ME_2)$.
\end{thm}

As a corollary we have the following branching rule for  $M^{n}_{k}.$ The assumption that ${\rm char}(F)\neq {\rm char}(\FF)$  is indispensable, it distinguishes  the cross-characteristic case from that  of defining characteristic. 

\medskip
\begin{thm}[Corollary 10.16\,~in~\cite{GDJ}]\label{IsoGrassTh1}
\quad If $0\leq {\rm char}(F)\neq {\rm char}(\FF)$ then
\begin{equation}\label{IsoGrass}M^{n}_{k}=M^{n-1}_{k-1}\,\oplus\, M^{n-1}_{k}\,\oplus\, \,M^{n-2}_{k-1}\,S_{n-1}\end{equation} as $FA_{n-1}$-modules.
\end{thm}

\subsection{Incidence maps}\label{Sec2.3}
 The containment relation among subspaces of $V=\FF^{n}$ yields an {\it incidence map\,} $\pa:M^{n}_{k}\to M^{n}_{k-1}$ for all $k.$ This map is defined on the basis $L^{n}_{k}$ of $M^{n}_{k}$ by setting $$\pa(x):=\sum\,\,y$$ where the sum runs over all co-dimension $1$ subspaces $y$ of $x,$ for $x\in L^{n}_{k}.$ If $i\geq 1$ is an integer then
 \begin{equation}\label{powerdel}\pa^{i}(x)=c(i)\sum\,\,y\end{equation} where the sum runs over all co-dimension $i$ subspaces $y$ of $x$ and where $c_{i}$ is a coefficient depending only on $i.$ It is easy to see that $c_{i}$ is the number of chains $y=x_{i}\subset x_{i-1}\subset\dots\subset x_{1}\subset x$ of length $i,$ and hence $$c(i)=[(i)!]_{q}=[i]_{q}\cdot[i-1]_{q}\cdots[1]_{q}\,\,.$$
 
Clearly $\p^{i}$ is an $FG_{n}$-homomorphism, and in particular,  the idempotents $E_{i}$ defined above commute with $\p^i.$ The following is useful for computations. If $x,\,y$ are subspaces of $V$ let $x\cdot y$ denote the subspace spanned by $x$ and $y.$ In particular, we write $v_{1}\cdot y$ for the subspace $\langle v_{1},\,y\rangle.$
This can be extended linearly to an associative and commutative product on $M^{n}:= \bigoplus M^{n}_{k}.$ 

Let $0\leq k\leq n$ and $f\in M^{n}_{k}.$ Then there are unique elements $f_{1}\in M^{n-1}_{k-1}$ and $\ell\in M^{n-1}_{k}$ so that $$f=v_{1}\cdot f_{1}\,+\,\ell\,.$$ This is  the {\it standard decomposition} of $f.$ We collect some useful facts. 

\begin{lem}\label{lem:231} (i)\, Let $f_{1}\in M^{n-1}_{k-1}.$ Then $(v_{1}\cdot f_{1})E_{1}=v_{1}\cdot f_{1}$ and $(v_{1}\cdot f_{1})E_{2}=0.$ \newline
(ii) \,Let $x=(a\big| x')\in L^{n}_{k}$ with $a\in(\FF^{k})^{T}$ and $x'\in L^{n-1}_{k}$ be a space not containing $v_{1}.$  Then $xE_{1}=q^{-k}\sum (b\big|x')$ where the sum runs over all $b\in(\FF^{k})^{T}.$ Furthermore, $xE_{2}=0$ unless $x'=v_{2}\cdot x''$ with $x''\in L^{n-2}_{k-1}.$
\end{lem}

\pf (i) The first part is a simple calculation, and for the second part we have $(v_{1}\cdot f_{1})E_{1}=v_{1}\cdot f_{1}$ so that $(v_{1}\cdot f_{1})E_{2}=(v_{1}\cdot f_{1})E_{1}E_{2}=0.$\, (ii) This is also a direct  calculation, alternatively see Theorem 10.2 \,in~\cite{GDJ} for a more general case.  \dne

\begin{lem}\label{lem:232} \,Let $f_{1}\in M^{n-1}_{k-1}$ and $i\geq1.$ Then $\pa^{i}(v_{1}\cdot f_{1})=v_{1}\cdot \pa^{i}(f_{1}) + q^{k-i}[i]_{q}\p^{i-1}(f_{1}E_{1}).$\end{lem}

\pf By linearity we can assume that $f_{1}=x'\in L^{n-1}_{k-1}$ and so  
$\pa^{i}(v_{1}\cdot x')$ represents the sum in $M^{n}_{k-i}$ of all $(k-i)$-dimensional subspaces of $v_{1}\cdot x'$ with coefficient $[(i)!]_{q},$ see \,(\ref{powerdel}). On the right hand side of the equation we have subspaces $y$ of $v_{1}\cdot x'$ of dimension $k-i,$ and $y$ has the correct coefficient $[(i)!]_{q}$ if it contains $v_{1}.$ Otherwise $y$ is a summand in the second term on the right hand side. In this case its coefficient is $q^{k-i}[i]_{q}[(i-1)!]_{q}\cdot q^{-(k-i)}$ see Lemma~\ref{lem:231}(ii). Now $q^{k-i}[i]_{q}[(i-1)!]_{q}\cdot q^{-(k-i)}=[(i)!]_{q}$ as required. \dne

\section{Partially Ordered Sets} 
 
 \vspace{-3mm}
 
 Let $(\CP,\leq )$ be a finite ranked poset with rank function ${\rm rk}\!:\,\CP\to\mathbb{N}\cup\{0\},$ and  assume that  $\min\{{\rm rk}(x)\,:\,x\in \CP\}=0.$ Let $G$ be the automorphism group of $\CP$ and let $F$ be a coefficient ring with $1.$ We describe a construction that associates to this data a family of homology modules which arise from the $FG$-action on certain permutation modules obtained from $\CP.$
 
For $k\in \mathbb{Z}$ denote the set of all $x\in {\cal P}$ with ${\rm rk}(x)=k$ by   ${\cal P}_{k}$ and let $M_{k}:=F{\cal P}_{k}$  be the $F$-vector space with basis ${\cal P}_{k}.$  In particular $M_{k}=0$  when $k<0$ or $k>\max\{{\rm rk}(x)\,:\,x\in \CP\}.$  Consider the linear {\it incidence map} $\partial\!: M_{k}\to M_{k-1}$ which is defined for $x\in  {\cal P}_{k}$ by $$\partial(x)=\sum y$$ where the sum runs over all $y\in{\cal P}_{k-1}$ with $y\leq x.$ Then \,$\partial$ is an $FG$-homomorphism on $M=\bigoplus_{k\in\mathbb{Z}}\,M_{k}$ and gives rise to the $FG$-sequence 
 \begin{equation}\label{SEQ00}
{\mathcal M}:\quad 0\stackrel{\,\,\,\partial}{\leftarrow}M_{0}\stackrel{\,\,\,\partial}{\leftarrow}
\ldots\stackrel{\,\,\,\partial}{\leftarrow}
M_{k-2}^{}\stackrel{\,\,\,\partial}{\leftarrow}M_{k-1}\stackrel{\,\,\,\partial}{\leftarrow}M_{k}\stackrel{\,\,\,\partial}{\leftarrow}M_{k+1}
\stackrel{\,\,\,\partial}{\longleftarrow}\ldots \stackrel{\,\,\,\partial}{\longleftarrow} 0\,.
\end{equation}
From this sequence  we obtain homological $FG$-sequences as follows. Note  that $M_{k}\neq 0$ for only finitely many $k,$ as $\CP$ is finite, and hence $\pa$ is nilpotent.  Hence let $m>1$ be an integer for which  $\pa^{m}=0$ as a map on $M.$ 
Fixing some  $0\leq k$ and  $0<i<m$ now consider the sequence  
\begin{equation}\label{SEQki}
{\mathcal M}^{}_{k,i}:\quad \ldots\stackrel{\,\,\,\partial^*}{\longleftarrow}
M_{k-m}^{}\stackrel{\,\,\,\partial^*}{\longleftarrow}
 M_{k-i}^{}\stackrel{\,\,\,\partial^*}{\longleftarrow}
M_{k}^{}
\stackrel{\,\,\,\partial^*}{\longleftarrow} M^{}_{k+m-i}
\stackrel{\,\,\,\partial^*}{\longleftarrow} M^{}_{k+m}
\stackrel{\,\,\,\partial^*}{\longleftarrow}\ldots\
\end{equation}
in which $\partial^{*}=\pa^{i}$ or $\pa^{m-i}$ stands for the appropriate power of $\partial$. Then ${\mathcal M}^{}_{k,i}$ is homological since $\partial^{*}\partial^{*}=\pa^{m}=0.$ We denote the
homology at $M_{k-i}^{}\leftarrow M_{k}^{}\leftarrow M^{}_{k+m-i}$
by
$$
H_{k,i}:=(\ker\,\partial^{i}\cap M_{k}^{})\,\,\Big/\,\,
                      \partial^{m-i}(M_{k+m-i}^{})\,\,.
$$
We call $H_{k,i}$ the {\it  incidence homology } of ${\cal P}$ over $F$ for parameters $k$ and $0<i<m.$ For example, when  $m=2$ then ${\mathcal M}={\mathcal M}_{k,1}$ for any $k\geq 0$ is a homological sequence in the usual sense.  This definition  of the homology depends on $m,$ and the choices for $m$ depend on both $F$ and $\CP.$  For instance, one  may take $m$ to be minimal with the property $\pa^{m}=0,$ but this is not the only case to consider.

The homology modules are related to each other. Denote  $I_{k,i}:=\partial^{m-i}(M_{k+m-i}^{})$ and $K_{k,i}:=(\ker\,\partial^{i}\cap M_{k}^{}).$ Then $\pa$ induces a linear map
\begin{equation}\label{hompa} \pa:H_{k,i}\to H_{k-1,i-1} \quad\text{by}\quad \pa(x+I_{k,i})=\pa(x)+I_{k-1,i-1}\,\,.
\end{equation}
In a similar way the identity map $\ic\!\!:M_{k}\to M_{k}$ induces a map $\ic\!\!:H_{k,i}\to H_{k,i+1}.$ If the $H_{k,i}$ are arranged as a grid with rows indexed by $0<i<m$ and columns indexed by $k$ then $\pa$ and $\ic$ connect the modules in  the array 
\begin{equation}\begin{array}{ccc} H_{k,i}&\quad&H_{k+1,i}\\&&\\
                                       & \nwarrow^{\pa}&\downarrow^{\ic}\\&&\\
                                         &&H_{k+1,i+1}
\end{array}\label{arrows}\end{equation}

of maps which carries essential information. (This feature is particular, it does not exist in ordinary homology.) It is used  to determine completely the homologies of finite projective spaces in Section~5.

We emphasize that this construction is entirely general, it applies to  an arbitrary finite ranked poset for any arbitrary coefficient ring with $1.$ 

\subsection{Branching Rules in Projective Space}

We now consider this homology for finite projective spaces.  Let  $q$ be a prime power, $\FF$ the field of $q$ elements and $n\geq 0.$ Then  the projective space $\CP=\CP(n,q)$ is the set  of all subspaces of $\FF^{n}$ ordered by inclusion, with rank function given by dimension.  Let $G_{n}:={\rm GL}(n,q).$ 
As coefficients we choose a field $F$ of characteristic $p_{}>0$ not dividing $q$ and let $$m=m(p_{},q)$$ be the characteristic of $\FF$ in $F$ discussed at the beginning of Section~\ref{Sec2}. 

\medskip
We  use the same notation for modules and incidence maps as above and as  in Sections~\ref{Sec2.2} and~\ref{Sec2.3}. It is clear from  
 \,(\ref{powerdel})\, that $\p^{m}=0$ on $M=\bigoplus_{k\in\mathbb{Z}}\,M_{k}$ independently of $n.$ In particular, if  $0\leq k$ and $0<i<m$ are given then we obtain the homological $FG$-sequence (or $FG$-chain complex)
\begin{equation}\label{SEQCAL}
{\mathcal M}^{n}_{k,i}:\quad \ldots\stackrel{\,\,\,\partial^*}{\longleftarrow}
M_{k-m}^{n}\stackrel{\,\,\,\partial^*}{\longleftarrow}
 M_{k-i}^{n}\stackrel{\,\,\,\partial^*}{\longleftarrow}
M_{k}^{n}
\stackrel{\,\,\,\partial^*}{\longleftarrow} M^{}_{k+m-i}
\stackrel{\,\,\,\partial^*}{\longleftarrow} M^{}_{k+m}
\stackrel{\,\,\,\partial^*}{\longleftarrow}\ldots\ \quad.
\end{equation}
The homology for parameters $k$ and $0<i<m$ is denoted by $$H^{n}_{k,i}:=\, K^{n}_{k,i}\,\Big/\,I^{n}_{k,i}.$$
It is convenient to set $K^{n}_{k,0}=0,$  $K^{n}_{k,m}=M^{n}_{k}$ and $I^{n}_{k,0}=0,$ $I^{n}_{k,m}=M^{n}_{k}.$ In particular, $H^{n}_{k,0}=0=H^{n}_{k,m}$ for all $k$ and $n.$

\bigskip 
It is appropriate  to discuss homology at the level of $FG$-sequences, rather than individual modules. (We prefer the term $`FG$-sequence' \, instead of   $`FG$-chain complex' \,as it  is simpler and does not interfere with the notion of a combinatorial complex.) 

Recall that if \,${\mathcal A}\!: \,\,... \stackrel{\,\,\,\alpha}{\longleftarrow}A_{k-1}\stackrel{\,\,\,\alpha}{\longleftarrow}
 A_{k}\stackrel{\,\,\,\alpha}{\longleftarrow}...$ \,\,and\, ${\mathcal B}\!: \,\,... \stackrel{\,\,\,\beta}{\longleftarrow}B_{k-1}\stackrel{\,\,\,\beta}{\longleftarrow}
 B_{k}\stackrel{\,\,\,\beta}{\longleftarrow}...$\, are finite $FG$-sequences then ${\cal A}\cong{\cal B}$ are isomorphic if and only if there is an $FG$-isomorphism  $\psi\!: \,A_{j}\to B_{j}$ for all $j$  such that $\psi\alpha=\beta\psi.$ (Assume that the first non-zero modules are $A_{0}$ and $B_{0}$ respectively.) Similarly, if ${\mathcal C}\!: \,\,... \stackrel{\,\,\,\gamma}{\longleftarrow}C_{k-1}\stackrel{\,\,\,\gamma}{\longleftarrow}
 C_{k}\stackrel{\,\,\,\gamma}{\longleftarrow}...$ is an $FH$-sequence for some subgroup $H\subseteq G$ and if  $A_{k}=C_{k}S$ is obtained by induction from $H$ to $G$ for each $k$ for some $S\subseteq G$ then ${\cal A}={\cal C}S$ is the $FG$-sequence induced from ${\cal C}$ with maps induced from $\gamma.$ 
 
\bigskip 
Let $S_{n-1}$ denote  a Singer cycle in $G_{n-1}={\rm GL}(n-1,q),$ of order $|S_{n-1}|= q^{n-1}-1,$ and let $M^{n-2}_{k-1}S_{n-1}$ be the module induced from $G_{n-2}$ to $G_{n-1}$ via Harish-Chandra induction. Correspondingly  let ${\cal M}^{n-2}_{k-1,i}S_{n-1}$ be the $FG_{n-1}$-sequence induced from   ${\cal M}^{n-2}_{k-1,i}\,.$
 
 \bigskip
\begin{thm} [Branching Rule]\label{thm:BRAM}
Let $\FF$ be the field of $q$ elements and let $F$ be a field of characteristic $p_{}>0$ not dividing $q.$  Suppose that $0 \leq  k \leq n$ with $1\leq n$ and $0<i<m=m(p_{},q)$ are integers.  Then
\[{\cal M}^n_{k,i}\,\, \cong \,\,{\cal M}^{n-1}_{k,i+1}\,\, \oplus\,\, {\cal M}^{n-1}_{k-1,i-1} \,\,\oplus\,\, {\cal M}^{n-2}_{k-1,i}\,S_{n-1}\]
as $FG_{n-1}$-sequences.   \end{thm} 

This result  is Theorem~\ref{IsoGrassTh1}\, at the level of $FG_{n-1}$-sequences; the proof is given in the next section. The following diagram may clarify the the situation. We show that there is an isomorphism $\vartheta$ of $FG_{n-1}$-sequences so  that all maps commute.

{\small \hspace{2.5cm}\begin{picture}(100,80)


\put(0,0){$M^{n}_{k-i}$}
\put(35,0){$M^{n-1}_{k-i}$}
\put(65,0){$M^{n-1}_{k-i-1}$}
\put(100,0){$M^{n-2}_{k-i-1}S_{n-1}$}

\put(0,33){$M^{n}_{k}$}
\put(35,33){$M^{n-1}_{k}$}
\put(65,33){$M^{n-1}_{k-1}$}
\put(100,33){$M^{n-2}_{k-1}S_{n-1}$}

\put(0,66){$M^{n}_{k-i+m}$}
\put(35,66){$M^{n-1}_{k-i+m}$}
\put(65,66){$M^{n-1}_{k-i+m-1}$}
\put(100,66){$M^{n-1}_{k-i+m-1}S_{n-1}$}


\put(17,1){\vector(1,0){12}}\put(19,3.5){\small $\vartheta_{k-i}$}\put(52,0.5){$\oplus$}\put(87,0.5){$\oplus$}
\put(17,34){\vector(1,0){12}}\put(19,36.5){\small $\vartheta_{k}$}\put(52,33.5){$\oplus$}\put(87,33.5){$\oplus$}
\put(17,67){\vector(1,0){12}}\put(19,69.5){\small $\vartheta_{k-i+m}$}\put(52,66.5){$\oplus$}\put(87,66.5){$\oplus$}

\put(2.5,24){\vector(0,-1){14}}\put(3.5,19){\small $\p^{i}$}
\put(105.5,24){\vector(0,-1){14}}\put(106.5,19){\small $\p^{i}$}

\put(48,24){\vector(1,-1){14}}\put(37.5,19){\small $\p^{i+1}$}
\put(60,24){\vector(-1,-1){14}}\put(59.5,19){\small $\p^{i-1}$}

\put(2.5,57){\vector(0,-1){14}}\put(3.5,52){\small $\p^{m-i}$}
\put(105.5,57){\vector(0,-1){14}}\put(106.5,52){\small $\p^{m-i}$}

\put(48,57){\vector(1,-1){14}}\put(33.5,52){\small $\p^{m-i+1}$}
\put(60,57){\vector(-1,-1){14}}\put(59.5,52){\small $\p^{m-i-1}$}

\end{picture}

\vspace{1cm}
\centerline{\sc Figure 1: Branching Rule}
}

Since ${\cal M}^{n}_{k,i}$ is homological we have corresponding $FG$-sequences for images and kernels, and this gives the sequence 
\begin{equation}\label{SEQHOM}
{\mathcal H}^{n}_{k,i}:\quad \ldots
\stackrel{\,\,\,\partial^*}{\longleftarrow}
H_{k-i,m-i}^{n}\stackrel{\,\,\,\partial^*}{\longleftarrow}
H_{k,i}^{n}
\stackrel{\,\,\,\partial^*}{\longleftarrow} H^{n}_{k+m-i,m-i}
\stackrel{\,\,\,\partial^*}{\longleftarrow} 
\ldots\ \quad.
\end{equation}
Using Theorem~\ref{thm:BRAM} and standard results from  homological algebra we have the following corollary: 

 \bigskip
\begin{thm} [Homology Decomposition]\label{thm:BRAH}
Let $\FF$ be the field of $q$ elements and let $F$ be a field of characteristic $p_{}>0$ not dividing $q.$  Suppose that $0 \leq  k \leq n$ and $0<i<m=m(p_{},q)$ are integers.  Then
\[{\cal H}^n_{k,i}\,\, \cong \,\,{\cal H}^{n-1}_{k,i+1}\,\, \oplus\,\, {\cal H}^{n-1}_{k-1,i-1} \,\,\oplus\,\, {\cal H}^{n-2}_{k-1,i}\,S_{n-1}\]
as $FG_{n-1}$-sequences.   \end{thm} 

This in turn completes the proof of  Theorem~\ref{bthm:HMD1}.
\medskip

\subsection{Proof of the Branching Rule}

The proof   requires two simple observations. First consider the idempotent $E_{1}$ introduced in Section~\ref{Sec2.2}. Observe that  $\p$ restricts to a map $\p\!:M^{n}_{k}E_{1}\to M^{n}_{k}E_{1}$ and hence we may consider the $FG_{n-1}$-sequence ${\cal M}^{n}_{k,i}E_{1}.$ We have 

\medskip
\begin{lem}\label{prop:BRAME1} Suppose that $0 \leq  k \leq n$ and $0<i<m$ are integers.  Then
\[{\cal M}^n_{k,i}E_{1}\,\, \cong \,\,{\cal M}^{n-1}_{k,i+1}\,\, \oplus\,\, {\cal M}^{n-1}_{k-1,i-1} \]
as $FG_{n-1}$-sequences.   \end{lem} 
 
\pf If $0\leq j\leq n$ is an integer and $f$ belongs to $M^{n}_{j}E_{1}$   consider its standard decomposition $f=v_{1}\cdot f_{1}+\ell$ as discussed in Section~\ref{Sec2.3}. Thus $f_{1}$ belongs to $M^{n-1}_{j-1}$ and none of the spaces appearing in $\ell$ contains $v_{1}.$ Hence $\ell$ is of the shape $$\ell=\sum\,\,\ell_{b,x'}\,\big(b\,\big|\,x')$$ 
where the sum runs over all $b\in (\FF^{k})^{\perp}$ and $x'\in L^{n-1}_{k},$ see Lemma~\ref{lem:231}. Hence let  $$\ell^{*}:=\sum\,\,\ell_{b,x'}\,\,x'\,\,\in M^{n-1}_{k}\,.$$ 
(Explicitly, $\ell^{*}$ is obtained by removing the first column of all matrices appearing in $\ell.)$ We define the map $\psi_{j}\!:M^{n}_{j}E_{1}\to M^{n-1}_{j}\,\oplus\,M^{n-1}_{j-1}$
for $j=k+m-i,\,\,k,\,\,k-i,...$ by
$$\begin{array}{lcccc}\psi_{k-i+m}(f)&=\,\,&\ell^{*}&\,\,+\,\,&q^{k}[m-i]_{q}f_{1}+\p(\ell^{*})\\\\

\psi_{k}(f)&=&\ell^{*}&+&q^{k-i}[i]_{q}f_{1}+\p(\ell^{*})\\\\

\psi_{k-i}(f)&=&\ell^{*}&+&q^{k-m}[m-i]_{q}f_{1}+\p(\ell^{*})\,,\,\,\,\, {\rm etc.}\,,
\end{array}$$
see  again Figure 1. It is clear that $\psi$ is an injective $FG_{n-1}$-homomorphism. By Theorems~\ref{thm:BR}\, and \ref{IsoGrassTh1}\, we have  $\dim M^n_{k}E_{1}=\dim M^{n-1}_{k}\,+\,\dim M^{n-1}_{k-1}$ and hence $\psi$ is an isomorphism. To prove  that $\psi$ is a homomorphism  of sequences it suffices to examine the case $j=k,$ the other indices are the same. To show that $(\pa^{i+1}\oplus\p^{i-1})(\psi_{k}(f))=\psi_{k-i}(\pa^{i}(f))$ we obtain 
\begin{eqnarray}
(\pa^{i+1}\oplus\p^{i-1})(\psi_{k}(f))&=&\p^{i-1}\big(q^{k-i}[i]_{q}f_{1}+\p(\ell^{*})\big)\,\,\,\,+\,\,\,\,\p^{i+1}(\ell^{*})\nonumber\\\nonumber\\
&=&q^{k-i}[i]_{q}\p^{i-1}(f_{1})+\p^{i}(\ell^{*})\,\,\,\,+\,\,\,\,\p^{i+1}(\ell^{*})\,.\nonumber
\end{eqnarray}

For the standard decomposition of $\p^{i}(f)$ we have 
$\p^{i}(v_{1}\cdot f_{1}+\ell)=v_{1}\cdot\pa^{i}(f_{1}) \,+\, q^{k-i}[i]_{q}\p^{i-1}(f_{1}E_{1})+\p^{i}(\ell)$ by Lemma~\ref{lem:232}\, and therefore
\begin{eqnarray}
\psi_{k-i}\big(\p^{i}(v_{1}\cdot f_{1}+\ell)\big)&=&q^{k-i}[i]_{q}\p^{i-1}(f_{1})+\p^{i}(\ell^{*})\,+\,
\big(q^{k-m}[m-i]_{q}+q^{k-i}[i]_{q}\big)\p^{i}(f_{1})\nonumber\\\nonumber&&+\,\,\,\,\,\,\p^{i+1}(\ell^{*}).
\end{eqnarray}
The result follows since $q^{k-m}[m-i]_{q}+q^{k-i}[i]_{q}=q^{k-m}(1+q+...+q^{m-1})=0$ in $F.$\dne 

\bigskip
Next we consider the idempotent $E_{2};$ in particular, assume that $F^{\times}$ contains an element of order $q_{0}$ where $q_{0}$ is the prime dividing $q.$   Here we have  the $FG_{n-2}$-sequence ${\cal M}^{n}_{k,i}E_{2},$ we claim that there is  the following isomorphism: 

\medskip
\begin{lem}\label{prop:BRAME2} Suppose that $0 \leq  k \leq n$ and $0<i<m$ are integers.  Then
\[{\cal M}^n_{k,i}E_{2}\,\, \cong \,\, {\cal M}^{n-2}_{k-1,i} \]
as $FG_{n-2}$-sequences.   \end{lem} 

\pf Let $0\leq j\leq n.$ Then the elements $(v_{2}\cdot x'')E_{2}$ with $x''\in L^{n-2}_{j-1}$ form a basis of $M^{n}_{j}E_{2},$ see Lemma~\ref{lem:232}. For $j=k+m-i,\,\,k,\,\,k-i,...$  define the map $\varphi_{j}\!:M^{n}_{j}E_{2}\to M^{n-2}_{j-1}$ by
$$\varphi_{j}((v_{2}\cdot x'')E_{2})=x''\,\,\in \,
M^{n-2}_{k-1}\,\,.$$
(Explicitly, $\varphi$ removes the first two columns in all matrices  appearing in  $f\in M^n_{k,i}E_{2}.)$ It is immediate that $\varphi$ is an $FG_{n-2}$-isomorphism, and similarly that $\varphi$ commutes with $\p^{i}$ and $\p^{m-i}$ as appropriate. \dne

\medskip
{\it Proof of Theorem~\ref{thm:BRAM}:}\,\, First assume that $F^{\times}$ contains an element of order $q_{0}.$ Let $S_{n-1}$ be a Singer cycle of $G_{n-1}.$ By the earlier comment on induced sequences and Lemma~\ref{prop:BRAME2}\,  we have an isomorphism ${\cal M}^n_{k,i}E_{2}S_{n-1}\, \cong \, {\cal M}^{n-2}_{k-1,i}S_{n-1}$ of $FG_{n-1}$-sequences which we also denote by $\varphi.$ Using Lemma~\ref{prop:BRAME1}\, we therefore have an isomorphism 
$$\vartheta:=(\psi,\varphi)\!:\,\,{\cal M}^n_{k,i}E_{1}\,\,\oplus\,\, {\cal M}^{n}_{k,i} E_{2}S_{n-1}
\,\,\,\cong \,\,\,\,{\cal M}^{n-1}_{k,i+1}\,\,\,\oplus\,\,\, {\cal M}^{n-1}_{k-1,i-1} \,\,\,\oplus\,\,\,{\cal M}^{n-2}_{k-1,i}S_{n-1}$$ of $FG_{n-1}$-sequences. Finally ${\cal M}^n_{k,i}={\cal M}^n_{k,i}E_{1}\,\,\oplus\,\, {\cal M}^{n}_{k,i}E_{2}S_{n-1}$ by Theorem~\ref{thm:BR}\, since $\p^{*}$ commutes with $E_{1}$ and $E_{2}.$ This proves the theorem when $F^{\times}$ contains an element of order $q_{0}.$ In the remaining case, extend $F$ to a field $\bar{F}\supset F$  which  contains an element of order $q_{0}$ and apply the result to this larger field. Now notice that the map $(\psi,\varphi)$ restricts back to  a map over $F.$   \dne

\section{Homology Modules $H^{n}_{k,i}$}

\vspace{-3mm}
We begin to analyze the incidence homology for finite projective spaces in detail.  As before $\FF={\rm GF}(q)$ and  $F$ is a field of characteristic $p_{}>0$ not dividing $q.$ Denote  $G_{n}={\rm GL}(n,q)$ and let $m=m(p_{},q)$ be the characteristic of $q$ in $F.$   

Let $n\geq 0.$ For any $k\leq n$ and $i$ with $0<i<m$ we have the  homological sequence 
\begin{equation} 
{\mathcal M}^{n}_{k,i}:\quad
0\stackrel{\,\,\,\gd^{*}}{\leftarrow}
...\stackrel{\,\,\,\gd^{*}}{\leftarrow}
M^{n}_{k-m}\stackrel{\,\,\,\gd^{*}}{\leftarrow}
M^{n}_{k-i}\stackrel{\,\,\,\gd^{*}}{\leftarrow}
M^{n}_{k}\stackrel{\,\,\,\gd^{*}}{\leftarrow}
M^{n}_{k+m-i}\stackrel{\,\,\,\gd^{*}}{\leftarrow}...\stackrel{\,\,\,\gd^{*}}{\leftarrow} 0\,\,
\end{equation}
in which  $\gd^{*}$ is the appropriate power of $\gd.$ Suppose that $n\geq m.$ Then for any choice of  $a<b$ in $\{0,..,m\!-1\}$ the sequence ${\mathcal M}^{n}_{b,b-a}$ contains the term $M^{n}_{a}\stackrel{\,\,\,\gd^{*}}{\leftarrow}M^{n}_{b},$ and all ${\mathcal M}^{n}_{k,i}$ are of this form. Hence there are ${m\choose 2}$ distinct sequences if $n\geq m.$ For instance, the index pairs $(k,i),$  $(k-i,m-i),$ and $(k\pm m,i)$ etc. all define the same sequence. We write $(k,i)\sim(k',i')$ if \,${\mathcal M}^{n}_{k,i}={\mathcal M}^{n}_{k',i'}.$ 

\subsection{Almost exact sequences  and Brauer characters}

\vspace{-3mm}
If ${\cal A}$ is  a homological sequence then its homology is {\it concentrated in a single position}, or\, ${\cal A}$ is {\it almost exact,} if all but at most one of the   homology modules in ${\cal A}$ vanish.  We show that ${\mathcal M}^{n}_{k,i}$ has this property for all $(k,i).$  One direction of the  following is Theorem~3.1\, in the paper~\cite{MSP} with Valery Mnukhin.

\medskip
\begin{thm}\label{thm:mtc} 
\, Let $0\leq k\leq n$ and  $0<i< m.$ Then $H^n_{k,i}\neq 0$ if and only if  $n <2k+m-i< n + m.$ \end{thm}

We call $(k,i)$  a {\it middle index} \,for  $n$ \,if (a) $0\leq k\leq n,$\, (b) $0<i< m$ and (c) $n <2k+m-i< n + m.$ If $(k,i)$ is a middle index then $M^{n}_{k-i}\leftarrow M^{n}_{k}$ is the {\it middle index} of the (unique) sequence  in which the two modules appear. It is easy to check that for any $k',\,i'$ the sequence  ${\mathcal M}^{n}_{k',i'}$ contains at most one middle index. But it may contain none. When  $m=2,$ for instance, then there is a unique middle index  $(k,1)$ for $n$ even, and $k=\frac n2$ in this case, while there is no middle index when $n$ is odd.

{\it Proof of Theorem~\ref{thm:mtc}}\,:\,  Assume that $0<i< m$ and $0\leq k\leq n.$ Writing out the inequalities observe that $(k,i)$ is a middle index for $n$ if and only if  $(k,i+1),$ $(k-1,i-1)$ and $(k-1,i)$ are middle indices for $n-1,$ $n-1$ and $n-2,$ respectively, unless $i\in \{1,m-1\}$ or $k=0.$ The result now follows from Theorem~\ref{bthm:HMD1}\, and induction on $n.$\dne 
 
\medskip
\begin{cor} \label{cor:mtc}  
\,  For all $(k,i)$ with $0<i<m$ the sequence ${\mathcal M}^{n}_{k,i}$ is almost exact. Furthermore,  ${\mathcal M}^{n}_{k,i}$ is  exact if and only $(k,i)\sim (k',i')$ implies that $(k',i')$ is not  a middle index for $n.$  
\end{cor}

This corollary implies that the incidence homology of $\CP(n,q)$ lies in  the Burnside ring of ${\rm GL}(n,q)$ over $F.$ For each $(k,i)$ the sequence  ${\mathcal M}^{n}_{k,i}$ is homological and therefore the Hopf-Lefschetz trace formula (see for instance Theorem 22.1, Chapter 2 in Munkres~\cite{Munk})\, says that  

\vspace{-0.5cm}
\begin{equation} \label{Hopf1}
\bigoplus_{t\in\mathbb{Z}}\,\,\big(\, H^{n}_{k+tm,i}-H^{n}_{k-i+tm,m-i}\,\big)\,\,=\,\,\bigoplus_{t\in\mathbb{Z}}\,\,\big(\, M^{n}_{k+tm}-M^{n}_{k-i+tm}\,\big)\,\,
\end{equation} 

\vspace{-0.5cm}
If $(k,i)$ is the middle index of ${\mathcal M}^{n}_{k,i}$ then $H^{n}_{k,i}$ is the only non-trivial  homology in  \,(\ref{Hopf1})\, and hence Corollary~\ref{cor:mtc} gives  

\begin{cor}[Trace Formula]\label{cor:trace} Let $(k,i)$ be a middle index. Then 
\begin{equation} \label{Hopf2}H^{n}_{k,i}\,=\,\bigoplus_{t\in\mathbb{Z}}\,\, M^{n}_{k+tm}\,
-\,\,\,\bigoplus_{t\in\mathbb{Z}}\,\,M^{n}_{k-i+tm}\end{equation} as $FG_{n}$-modules in the Burnside ring. 
\end{cor}

Considering characters, if $H$ is any $FG_{n}$-module, let $\chi(g,\,H)$ denote the Brauer character of $G_{n}$ on $H.$ In particular, the permutation character
$$\chi(g,\,M^{n}_{k})=\pi_{k}(g)$$
is the number of $k$-dimensional spaces of $V$ that are stabilized by $g\in G_{n}.$ Hence Corollary~\ref{cor:trace} \, yields

\medskip
\begin{thm} 
\label{thm:Brauer} 
\,\,\, Let $(k,i)$ be a middle index for $n.$ Then the Brauer character of ${\rm GL}(n,q)$  on $H^n_{k,i}$ is  $$\chi(g,\,H^n_{k,i})\,\,=\,\,\sum_{t\in\mathbb{Z}}\,\,\,\,\pi_{k+tm}(g) - \pi_{k-i+tm}(g)\,.$$ \end{thm}

This key fact is already mentioned in ~\cite{MSP},
 it also provides the dimensions of the homology modules in the next section.

\subsection{Betti numbers}

\vspace{-3mm}
By Theorem~\ref{thm:mtc} we have  $H^n_{k,i}\neq 0$ if and only if  $(k,i)$ is a middle index, and in that case  let   $$\gb^{n}_{k,i}:=\dim H^{n}_{k,i}$$
be the Betti number of $H^n_{k,i}.$  As before it will be useful to set $\gb^{n}_{k,0}=0=\gb^{n}_{k,m}.$

\medskip\medskip
\begin{thm} \label{thm:Betti}
\quad Let $(k,i)$ be a middle index for $n.$ Then 

\vspace{-6mm}
\begin{enumerate}[\!\!\!(i)\,\,]
\item $\gb^n_{k,i}=\gb^{n-1}_{k,i+1}\,\,+\,\,\gb^{n-1}_{k-1,i-1}\,\,+\,\, \gb^{n-2}_{k-1,i}(q^{n-1}-1)\,;$  \vspace{-2mm}
\item $\gb^n_{k,i}\,=\,\sum_{t\in\mathbb{Z}}\,\,{n\choose k+tm}_{\!q}-{n\choose k-i+tm}_{\!q}.$ In particular, $\gb^n_{k,i}$ is the Euler characteristic of  ${\mathcal M}^{n}_{k,i};$\vspace{-2mm}
 \item
 \hspace{-1mm}{(Duality:)}\, For all $0\leq \ell\leq n$ and $0<j<m$ we have $\gb^n_{\ell,j}=\gb^n_{n-\ell,m-j}.$\vspace{-2mm}
\item
Let $\epsilon=2$ when  $n$ is odd and  $\epsilon=1$ otherwise.  Then as a polynomial in $q$ we have $$\gb^n_{k,i}\,=(q^{n-1}-1)(q^{n-3}-1)(q^{n-5}-1)\cdots(q^{\epsilon}-1)+ f(q)$$ where $f(q)$ is a polynomial of  degree $< (n-1)+(n-3)+(n-5)+\dots+\epsilon.$
\end{enumerate}
\end{thm}

{\sc Remark:} We say that  the index pair $(n-\ell,m-j)$ is  {\it dual}\, to $(\ell,j).$ The equation $\gb^n_{\ell,j}=\gb^n_{n-\ell,m-j}$ in (iii) extends to an isomorphism $H^n_{\ell,j}\,\cong\,H^n_{n-\ell,m-j}.$ This is shown in Section 5 which contains many  additional inequalities for Betti numbers.

\medskip
\pf The first part (i) follows from Theorem~\ref{bthm:HMD1}\, and (ii) follows from \,(\ref{Hopf2}). To show (iii) first notice that $(\ell,j)$ is a middle index if and only if the same is true about its dual. Hence we assume that $(\ell,j)$ is a middle index and proceed by induction. The equality $\gb^n_{\ell,j}=\gb^n_{n-\ell,m-j}$ is easily checked for $n\leq 2.$ Using (i) we have 
$$\gb^n_{\ell,j}=\gb^{n-1}_{\ell,j+1}+\gb^{n-1}_{\ell-1,j-1}+ \gb^{n-2}_{\ell-1,j}(q^{n-1}-1)$$ and 
$$\gb^n_{n-\ell,m-j}=\gb^{n-1}_{n-1-\ell,m-j-1}+\gb^{n-1}_{n-1-\ell+1,m-j+1}+ \gb^{n-2}_{n-2-\ell+1,m-j}(q^{n-1}-1)\,. $$ Observe that 
$(\ell-1,j)$ is dual to  $(n-2-\ell+1,m-j)$  relative to $n-2$ and the other two pairs  are dual to each other relative to $n-1.$ The result now follows by induction.
The last part follows by induction on $n$ from the first part. \dne

For small $m(p_{},q)$ one can analyze the Betti numbers a little further. 

\medskip
\begin{cor} \label{cor:HomDim}
\, (a) Suppose that $m(p_{},q)=2.$ If $n$ is odd then ${\mathcal M}^{n}_{k,1}$ is exact. If $n=2k$ is even then $\gb^{n}:=\gb^{n}_{k,1}$ depends on $n$ only and $$\gb^{n}=(q^{n-1}-1)(q^{n-3}-1)(q^{n-5}-1)\cdots(q-1).$$ 
(b) Suppose that $m(p_{},q)=3.$ For given $n$ either ${\mathcal M}^{n}_{k,i}$ is exact or $(k,i)$ is one of two middle indices $(k,1)\neq (k',2)$  when  $\gb^{n}:=\gb^{n}_{k,1}=\gb^{n}_{k',2}$ depends on $n$ only. Furthermore, 
$$\gb^{n}=\gb^{n-1}+\gb^{n-2}(q^{n-1}-1)$$
with initial values $\gb^{0}=\gb^{1}=1.$ \end{cor}

\pf We leave this to the reader. \dne

\section{Module Structure  and Duality}

We determine the incidence homology up to isomorphism, in terms of the standard ${\rm GL}(n,q)$-irreducibles over $F.$ In addition we show that there are dualities between homology modules for  certain index pairs.  

\subsection{Composition Factors}
If $(k,i)$ is a middle index for $n$ then $n<2k+m-i<n+m.$ We say that  $(k,i)$ is a {\it maximal middle index } if $2k-i=n-1.$ The following is checked easily.

\medskip
\begin{lem}\label{inductNEW} Let $n\geq 2$ and let $(k,i)$ be a maximal middle index for $n.$ Then\newline
{\rm (i)} \,\,\, $(k-1,i)$ is a maximal middle index for $n-2$ unless $k=0;$\\ 
{\rm (ii)} \,\, $(k,i+1)$ is a maximal middle index for $n-1$ unless $i+1=m;$\\ 
{\rm (iii)} \, $ (k-1,i-1)$ is a maximal middle index for $n-1$ unless $i-1=0$ or $k=0.$
\end{lem}

If  $\lambda=(\lambda_{1},\lambda_{2})$ is a composition of $n$ denote the Specht module for $\lambda$ by $S^{\lambda}$ and let $D^{\lambda}:=S^{\lambda}/S^{\lambda}\cap (S^{\lambda})^{\perp}$ be its head.  Then $D^{\lambda}= D^{(\lambda_{2},\lambda_{1})}$ is the usual standard irreducible $FG_{n}$-modules indexed by $\lambda.$ If $\mu=(\mu_{1},\mu_{2})$ is another  composition of $n$ then $\lambda$ majorizes $\mu,$ denoted $\lambda  \geq \mu,$ if a largest part of $\lambda$ is bigger or equal to a largest part of $\mu.$ We write $\lambda>\mu $ if $\lambda\geq \mu $ and $\lambda\neq \mu,$  see James' book~\cite{GDJ}.

\medskip
\begin{prop} \label{prop:IsoFactor}\, (i) \,Let  $(k,i)$ be a middle index.  Then all composition factors of $H^{n}_{k,i}$ have the form $D^{\lambda}$ with $\lambda\geq (n-k,k).$ Furthermore,  $D^{(n-k,k)}$ has multiplicity $1$ in $H^{n}_{k,i}.$

\vspace{-3mm}
(ii) \, If $(k,i)$ is a maximal  middle index for $n$  then $H^{n}_{k,i}\cong D^{(n-k,k)}.$ \end{prop}

\pf  (i) By Corollary 16.3\, in~\cite{GDJ}\, the composition factors of $M^{n}_{k}$  are of the form $D^{\lambda}$ with $\lambda\geq (n-k,k)$ and hence  the same is true for $H^{n}_{k,i}.$ Furthermore, $D^{(n-k,k)}$ has multiplicity $1$ in $M^{n}_{k}.$ The remainder is easy to verify when $n\leq 2$ and when $k=0.$ In the latter case we have $H^{n}_{k,i}=F=D^{(n,0)}.$ If $n\geq 2$ and if $(k,i)$ is a middle index of $n$ then $(k-1,i)$ is a middle index of $n-2.$ By induction and the homology decomposition we may assume that $D^{(n-1-k,k-1)}$ is a composition factor of $H^{n-2}_{k-1,i}.$

Let $D^{\lambda_{t}},$ $D^{\lambda_{t-1}},...$ $D^{\lambda_{1}}$ be the composition factors of $H^{n}_{k,i},$ ordered so that $\lambda_{t}>\lambda_{t-1}>...>\lambda_{1}\geq (n-k,k).$ For $\lambda=(a,b)$ with $a,b>0$ let $\lambda''$ be the composition $(a-1,b-1)$ of $n-2.$ It follows from Theorem 16.9\, of \cite{GDJ}\, that the composition factors of 
 $H^{n-2}_{k-1,i}=H^{n}_{k,i}E_{2}$ are of the shape $D^{\lambda_{j}''}$ with $\lambda_{j}''\geq (n-k,k)''=(n-1-k,k-1).$ As $D^{(n-1-k,k-1)}$ is a composition factor of $H^{n-2}_{k-1,i}$ we have $\lambda_{1}= (n-k,k).$  By Corollary 16.3\, 
 the multiplicity of $D^{(n-k,k)}$ in $H^{n}_{k,i}$ is $1.$ 
 
\smallskip
(ii) This  is easy to verify when $n\leq 2$ or when $k=0.$ So suppose that $n\geq 2$ and  $k\neq0.$ Using Lemma~\ref{inductNEW}\, we may apply induction to the terms on the right hand side of the decomposition $$H^n_{k,i} \,\cong\, H^{n-1}_{k,i+1} \,\oplus \,  H^{n-1}_{k-1,i-1} \,\oplus\, H^{n-2}_{k-1,i}\,S_{n-1}.$$ By the first part of the theorem we have  $H^{n-2}_{k-1,i}=D^{(n-k,k)''}$ since $k>0.$ Furthermore, we have $H^{n-1}_{k,i+1}=D^{(n-k,k)^{a}}$(unless $H^{n-1}_{k,i+1}=0$ when $i=m-1)$\, and $H^{n-1}_{k-1,i-1}=D^{(n-k,k)^{b}}$ (unless  $H^{n-1}_{k-1,i-1}=0$ when $i=1)$ where  $(n-k,k)^{a}:=(n-1-k,k)$ and   $(n-k,k)^{b}:=(n-k,k-1)$ are the corresponding compositions of $n-1.$ By the first part we have that $D^{(n-k,k)}$ has multiplicity one in $H^{n}_{k,i}.$
Suppose that also $D^{\lambda}$ with $\lambda >(n-k,k)$ appeared  as a composition factor in $H^{n}_{k,i}.$ Now use Theorem 16.9~of~\cite{GDJ}\, (or direct computation) to show that at least one of the modules $H^{n-1}_{k,i+1},$ $H^{n-1}_{k-1,i-1},$ $H^{n-2}_{k-1,i}$ has more than one composition factor, a contradiction. \dne

\subsection{Embeddings}

Let $0\leq k\leq n$ and $m=m(p_{},q)$ be given. Then  $(k,i)$ is a middle index  for  some $i$ if and only if $\frac12(n-m+1)<k<\frac12(n+m-1)$ and $0\leq k\leq n.$ Let $\ell$ and $r$ be the minimal and maximal choices for $k$ satisfying this constraint, respectively. It is useful to imagine the middle indices as nodes in an array of $m-1$ rows and $m-1$ or $m-2$ columns, depending on the parity of $n-m,$ if $n\geq m-2.$  (There will be fewer columns if $n<m-2$ due to the constraint $0\leq k\leq n.)$ Placing the non-zero homologies into this array we obtain the table 

{\small $${\cal H}^{n}\!:\,\quad \,\,\,\begin{array}{cccccccc} \medskip
H^{n}_{\ell,1} & H^{n}_{\ell+1,1} & H^{n}_{\ell+2,1} & \cdots && 0 & 0 & 0 \\\medskip
0 & H^{n}_{\ell+1,2} & H^{n}_{\ell+2,2} & \cdots && 0 & 0 & 0 \\ \medskip
0 & H^{n}_{\ell+1,3} & H^{n}_{\ell+2,3} & \cdots && 0 & 0 & 0 \\ \medskip
0 & 0 & H^{n}_{\ell+2,4} & \cdots && 0 & 0 & 0 \\ \medskip
0 & 0 & H^{n}_{\ell+2,5} & \cdots && 0 & 0 & 0 \\ \medskip
\vdots & \vdots & \vdots & & & \vdots & \vdots & \vdots \\ \medskip
0 & 0 & 0 && \cdots & H^{n}_{r-2,m-5}& 0 & 0 \\ \medskip
0 & 0 & 0 && \cdots & H^{n}_{r-2,m-4}& 0 & 0 \\ \medskip
0 & 0 & 0 && \cdots & H^{n}_{r-2,m-3} & H^{n}_{r-1,m-3} & 0 \\ \medskip
0 & 0 & 0 && \cdots & H^{n}_{r-2,m-2} & H^{n}_{r-1,m-2} & 0 \\ \medskip
0 & 0 & 0 && \cdots & H^{n}_{r-2,m-1} & H^{n}_{r-1,m-1} & H^{n}_{r,m-1}
\end{array}  
$$

\medskip
\centerline{Figure 2: ${\cal H}^{n}$-Table}
}

when $n-m$ is even, and a similar table with two non-zero entries in the first and last column, when $n-m$  is odd.  (Similarly, remove an equal number of columns on the left and the right of the array if $n<m-2.)$ The corner entries $H^{n}_{r,m-1},$ $H^{n}_{r-1,m-3},$  $H^{n}_{r-2,m-5}...$ on the right correspond to maximal middle indices. As discussed in Section 3\, the incidence map $\pa\!:\, M^{n}_{k}\to M^{n}_{k-1}$ induces  a map $\pa\!:\,H^{n}_{k,i}\to H^{n}_{k-1,i-1}.$ This map is  a NW-arrow in the table.  The identity map $M^{n}_{k}\to M^{n}_{k}$ induces  a map $\ic\!:\,H^{n}_{k,i}\to H^{n}_{k,i+1}, $ this  represents a S-arrow in the ${\cal H}^{n}$-table.  Compare  to the general situation described in \,(\ref{arrows}).

\medskip
\begin{lem} \label{maps}\, (i) Let $0\leq k\leq n$ and $1\leq t<i<m.$ Assume that  $2k-t\geq n.$ Then $\pa^{t}\!: H^{n}_{k,i}\to H^{n}_{k-t,i-t}$ is an injection. \newline
(ii) \, Let $0\leq k\leq n$ and $1<i<j<m.$ Assume  that  $2k+m-i-j\geq n.$ Then $\ic^{j-i}\!: H^{n}_{k,i}\to H^{n}_{k,j}$ is an injection. 
\end{lem}

Writing $2k-t=k+(k-t)$ the condition in (i) says that the arrow $H^{n}_{k,i}\to H^{n}_{k-t,i-t}$ is balanced on or towards the right of the centre of the array. Similarly, the assumption in  (ii) is a balancing condition on the diagonal of the array. 

\pf (i)  
We can assume that $n<2k+m-i<n+m$ since otherwise $H^{n}_{k,i}=0$ by Theorem~\ref{thm:mtc}, and in this case the assertion is true. Now $(k,t)$ is a middle index if and only if $2k-t<n.$ From the assumption and Theorem~\ref{thm:mtc}\, we conclude that $H^{n}_{k,t}=0.$ For injectivity we show that if  $\pa^{t}(x)$ belongs to  $\pa^{m-i+t}(M^{n}_{k+m-i})$ then  $x$ belongs to $\pa^{m-i}(M^{n}_{k+m-i}).$ So suppose that $\pa^{t}(x)=\pa^{m-i+t}(w)$ for some $w$ in $M^{n}_{k+m-i}.$ Then  $\pa^{t}\big(x-\pa^{m-i}(w)\big)=0$ and since $K^{n}_{k,t}=I^{n}_{k,t}$ there is  some $u$ in $M_{k+m-t}$ so that  $x-\pa^{m-i}(w)=\pa^{m-t}(u).$ Hence $x=\pa^{m-i}\big(w+\pa^{i-t}(u)\big).$ 

(ii) 
Again assume that $n<2k+m-i<n+m.$ Now $(k+m-j,m-j+i)$ is a middle index if and only if  $2k+m-i-j<n.$ From the assumption and Theorem~\ref{thm:mtc}\, we have $H^{n}_{k+m-j,m-j+i}=0.$ So suppose that $x=\pa^{m-j}(w)$ for some $w\in M^{n}_{k+m-j}$ and $\pa^{i}(x)=0.$ Then $\pa^{m-j+i}(w)=0$ and as $H^{n}_{k+m-j,m-j+i}=0$ we have $w=\pa^{j-i}(u)$ for some $u.$ Therefore $x=\pa^{m-j}(\pa^{j-i}(z))=\pa^{m-i}(u),$ and hence  $\ic^{j-i}\!: H^{n}_{k,i}\to H^{n}_{k,j}$ is injective.\dne

\medskip
\begin{thm}[Duality]\label{Dual} \,
Let $\FF={\rm GF}(q)$  and let $F$ be a field of characteristic $p_{}>0$ not dividing $q.$  Suppose that $0 \leq  k \leq n$ and $0<i<m=m(p_{},q).$  Then \newline 
(i) \,\, $H^n_{k,i}\, \cong \,H^{n}_{n-k,m-i},$ and \newline
(ii) \, $H^n_{k,i}\, \cong \,H^{n}_{k,j}$ where $j=2k-n+m-i$\newline
as $FG_{n}$-modules.
 \end{thm}

This is Theorem~\ref{bthm:Dual}\, in the Introduction. It endows the ${\cal H}^{n}$-array with a $C_{2}\times C_{2}$ symmetry which we will use to complete our analysis. 

{\sc Example:} \, To illustrate the  two dualities  consider the case $n=10,$$m=5,$ $k=4$ and $i=2.$ Then $j=1$ and the relevant index pairs are  $(4,2),$ $(6,3),$ $(4,1)$ and $(6,4).$ From the trace formula Corollary~\ref{cor:trace} we can express the $H^{n}_{k,i}$ as

\vspace{-5mm}
\begin{enumerate}[\qquad\qquad\qquad(a)\,\,]
\item  $H^{10}_{4,2}=M^{10}_{4}-M^{10}_{2}+M^{10}_{9}-M^{10}_{7},$\vspace{-2mm}
\item  $H^{10}_{6,3}= M^{10}_{6}-M^{10}_{3}+M^{10}_{1}-M^{10}_{8},$\vspace{-2mm}
\item  $H^{10}_{4,1}= M^{10}_{4}-M^{10}_{3}+M^{10}_{9}-M^{10}_{8},$\vspace{-2mm}
\item  $H^{10}_{6,4}= M^{10}_{6}-M^{10}_{2}+M^{10}_{1}-M^{10}_{7}.$\vspace{-2mm}
\end{enumerate}

\vspace{-3mm}
At the level of permutation sets $L^{n}_{k}$ is not permutation equivalent to $L^{n}_{n-k},$ and so the isomorphism does not hold for permutation sets. However, in the case of cross characteristics we do have  $M^{n}_{k}\cong M^{n}_{n-k}$ at the level of permutation modules, see Theorem~14.3~in~\cite{GDJ}. From this it follows that the four modules are indeed isomorphic.

\medskip 
{\it Proof of Theorem~\ref{Dual}:\quad } For the index pairs $(k,i)$ and $(\ell,j)$ we write $(k,i)\ra(\ell,j)$ or $(k,i)\lra(\ell,j)$ provided that there  is an $FG_{n}$-monomorphism, or  $FG_{n}$-isomorphism $H^n_{k,i}\to H^{n}_{\ell,j},$ respectively. For the first part of the theorem note that $\dim(H^n_{k,i})=\dim(H^{n}_{n-k,m-i})$  by Theorem~\ref{thm:Betti}(iii). In particular, $(k,i)$ is a middle index if and only if $(n-k,m-i)$ is a middle index. Hence it suffices to show that $(k,i)\ra(n-k,m-i)$ for all $k\geq \frac 12 n$ and all $0<i<m.$

{\sc When $n$ is even:} \,\, Here ${\cal H}^{n}$ has a middle column indexed by $\frac n2.$ Consider the middle index $(k,i)$ where $k=\frac 12 n +a$
with $0\leq a<\frac 12 i\leq \frac 12 (m-1)$ since $k<\frac 12 (n+m-1).$ First suppose that  $2i\leq m+2a.$ Then $\ic^{m+2a-2i}\!: (k,i)\ra (k,m-i+2a )$ and $\pa^{2a}\!:(k,m-i+2a )\ra (k-2a,m-i)$ by  Lemma~\ref{maps}. Hence $(k,i)\ra (n-k,m-i).$ Next suppose that $a>0$ and $2i> m+2a.$ Here we have  $\pa^{2a}\!:(\frac 12 n+a,i )\ra (\frac 12 n-a,i-2a)$ by the lemma. Now observe that $(\frac 12 n-a,i-2a)\lra (\frac 12 n+a,m-i+2a)$ by the first part of the proof, since $2(m-i+2a)\leq m+2a.$ Next we have $\pa^{2a}\!:(\frac 12 n+a,m-i+2a)\ra (\frac 12 n-a,m-i+2a-2a).$ Hence together $(\frac 12 n+a,i )\ra(\frac 12 n-a,m-i)$ which proves the result for $n$ even.  

{\sc When $n$ is odd:} \,\, The argument is almost the same. Consider the middle index $(k,i)$ where $k=\frac 12 n +a$ with $a=\frac 12,\,1+\frac 12,...,<\frac 12 i\leq \frac 12 (m-1).$  First suppose that $2i\leq m+2a.$ Then $\pa^{2a}\!:(k,i)=(\frac 12 n+a,i)\ra (\frac 12 n-a,i-2a)$ and   $\ic^{m-2i+2a}\!:(\frac 12 n-a,i-2a)\ra (\frac 12 n-a,m-i)$ by Lemma~\ref{maps}. Hence $(k,i)\ra (n-k,m-i)$ in this case. Next suppose that $2i> m+2a$ so that $\pa^{2a}\!:(k,i)=(\frac 12 n+a,i)\ra (\frac 12 n-a,i-2a).$ Now $(\frac 12 n-a,i-2a)\lra (\frac 12 n+a,m-i+2a)$ by the first part of the proof, as $2(m-i+2a)\leq m+2a.$ This gives $\pa^{2a}\!:(\frac 12 n+a,m-i+2a)\ra (\frac 12 n-a,m-i+2a-2a).$ Hence  $(\frac 12 n+a,i )\ra(\frac 12 n-a,m-i),$ and this completes the proof of the first part.

\smallskip
To prove the second part note that $(k,i)$ is a middle index if and only if $(k,j)$ with $0<j=2k-n+m-i<m$ is a middle index. Assume therefore  that $(k,i)$ is a middle index so that $0<2k-n+m-i<m$ and without loss $i<j.$ First suppose $k\geq \frac 12 n,$ with $k= \frac 12 n+a.$ Then $\ic^{2k-n+m}\!:(k,i)\ra(k,2k-n+m-i)=(k,j)$ by Lemma~\ref{maps}. Conversely, we have $\pa^{2a}\!: (k,2k-n+m-i)\ra(k-2a,2k-2a-n+m-i)$ by the lemma and $(k-2a,2k-2a-n+m-i)\lra(n-k+2a,m-2k+2a+n-m+i)$ by the first part of the theorem. But $(n-k+2a,m-2k+2a+n-m+i)=(k,i).$ Hence $(k,j)\ra (k,i).$ The case $k\leq \frac 12 n$ is very similar.  \dne

\medskip
By the last theorem all homologies are determined by the modules $H^{n}_{k,i}$ in the `triangle' where $2k=n+2a$ with $a\geq 0$  and $2i\leq  m +2a.$ It remains to examine these terms. 

\medskip
\begin{lem} \label{Hsum}Let $(k,i)$ be a maximal middle index for $n\geq 0,$ and write $k=\frac 12 n+a$ with $a\geq 0.$ Suppose that $j$ is an integer with $i<j<m$ and $2j\leq m+2a.$  Then $\gb^{n}_{k,j}=\gb^{n}_{k,i}+ \gb^{n}_{k+1,j+1}.$
\end{lem}

\pf Note, our usual convention applies, we put $\gb^{n}_{k+1,j+1}=0$ if $(k+1,j+1)$ is not a middle index.  The equation holds for $j=m-1$ by Theorem~\ref{Dual}(ii), and similarly for $a\geq \frac 12(m-3)$ in which case $k$ indexes the last column in the ${\cal H}^{n}$-array. The statement is also true for $0\leq n\leq 1.$ The result now follows by induction and Theorem~\ref{thm:Betti}(i).  \dne

\medskip

\begin{thm}\label{SumDual} \,
Let $\FF={\rm GF}(q)$  and let $F$ be a field of characteristic $p_{}>0$ not dividing $q.$ Let $(k,i)$ be a maximal middle index for $n\geq 0,$ and write $k=\frac 12 n+a$ with $a\geq 0.$ Suppose that $j$ is an integer with $i<j<m$ and $2j\leq m+2a.$
Then $H^{n}_{k,j}\,\cong\, H^{n}_{k,i}\,\,\oplus\,\, H^{n}_{k+1,j+1}$ as $FG_{n}$-modules.  
\end{thm}

\pf By Lemma~\ref{maps} we have  $\ic^{j-i}\!:(k,i)\ra(k,j)$ and $\pa\!:(k+1,j+1)\ra (k,j).$ As $(k,i)$ is a maximal middle index it follows  from Proposition~\ref{prop:IsoFactor}\, that $H^{n}_{k,i}=D^{(n-k,k)}$ is irreducible and  not a composition factor of $H^{n}_{k+1,j+1}.$  Therefore the sum $\ic^{j-i}(H^{n}_{k,i})+\pa(H^{n}_{k+1,j+1})\subseteq H^{n}_{k,j}$ is direct and the result follows from Lemma~\ref{Hsum}.\dne

\medskip

\begin{thm}\label{SumDual2} \,
Let $\FF={\rm GF}(q)$  and let $F$ be a field of characteristic $p_{}>0$ not dividing $q.$ Let $(k,j)$ be a middle index for $n\geq 0.$ Assume that  $k\geq\frac 12 n$  and  $j\leq \frac 12 (m-n) +k.$  Put $\ell=n-k+j-1.$
Then $H^{n}_{k,j}= \bigoplus_{t=k}^{\ell}\, D^{(n-t,t)}.$
\end{thm}

When $j>k+1$ then for some terms in $\bigoplus_{t=k}^{\ell}\, D^{(n-t,t)}$ we have expressions $(n-t,t)$ with $n-t<0.$ Here evidently $D^{(n-t,t)}=0.$

\pf First notice that $\ell\geq k$ since $(k,j)$ is a middle index, with equality if and only if $(k,j)$ is a maximal middle index. Next let  $k$ be maximal with $k<\frac 12 (n+m-1),$ that is, $k$  indexes the right-most column of ${\cal H}^{n}.$ 
If $\frac 12 (n+m-1)\leq n$ then the constraint $j\leq \frac 12 (m-n) +k$ implies that $(k,j)$ is a maximal middle index, and hence $k=\ell.$ From Proposition~\ref{prop:IsoFactor}\, it follows that  $H^{n}_{k,j}=D^{(n-k,k)}$ as required. If $k=n$ then $H^{n}_{k,j}=F=D^{(0,n)}=\bigoplus_{t=k}^{\ell}\, D^{(n-t,t)}$ having in mind the comment about expressions $(n-t,t)$ with $n-t<0.$ 
So we may apply induction and suppose that the theorem holds for all values $>k.$ Now apply Theorem~\ref{SumDual}\, and Propostion~\ref{prop:IsoFactor}. \,\dne

{\it Proof of Theorem~\ref{bthm:Comp}:} \,The first part is Theorem~\ref{thm:mtc} and so we turn to the second part. Let $(k,i)$ be a middle index. (a) If $k\geq \frac 12 n$ and  $i\leq \frac 12 (m-n) +k$ set $$T_{k,i}:=\{\,t\!:\,k\leq t\leq n-k+i-1\}.$$ In this case the statement of Theorem~\ref{bthm:Comp}\, is Theorem~\ref{SumDual2}.    Considering the remaining three possibilities suppose\, (b) that $k\geq  \frac 12 n$ and  $i>\frac 12 (m-n) +k.$ Here let $i'=2k-n+m-i.$ Then $H^{n}_{k,i}\,\cong\,H^{n}_{k,i'}$ by Theorem~\ref{Dual}(ii)\, and as $(k,i')$ is of type (a) the result follows if we set  $$T_{k,i}=\{\,t\!:\,k\leq t\leq k+m-i-1\}.$$ Similarly, using Theorem~\ref{Dual}(i)\, the remaining two cases can be reduced to (a) or (b) if we define 

\vspace{-4mm}
$$\begin{array}{rlll} T_{k,i}&=&\{\,t\!:\,n-k\leq t\leq n-k+i-1\}\,&\text{if $k<\frac 12 n$ \,and\, $i\leq \frac 12 (m-n)+k,$ and }\\\smallskip
T_{k,i}&=&\{\,t\!:\,n-k\leq t\leq k+m-i-1\}\,&\text{if $k<\frac 12 n$ \,and\, $i>\frac 12 (m-n)+k.$}\\
\end{array}$$
This completes the proof. In practical terms the result says that  the ${\cal H}^{n}$-array is bordered by $D^{(n-k,k)}$ on both ends of column $k$ with the remainder filled in using the simple rule of Theorem~\ref{SumDual2}\, and its  mirror versions. \dne

\end{document}